\documentclass[11pt]{article}
\usepackage{float}
\usepackage{graphicx}
\usepackage{enumerate}
\usepackage{amsmath}
\usepackage{amsfonts}
\usepackage{amssymb}
\usepackage{amsfonts}
\usepackage{amssymb,amsmath,amsthm,epsfig,graphicx,subfigure,hyperref}
\usepackage[margin=1.3in]{geometry}
\usepackage{caption2}
\usepackage{color}
\usepackage{enumerate}
\usepackage{booktabs}
\usepackage[square, sort, numbers]{natbib}
\usepackage{supertabular}
\usepackage{blkarray}
\usepackage{authblk}
\RequirePackage{lineno}

\def\blackbox{\hfill {\vrule height3pt width4pt depth2pt}}
\def\Box{\hfill \framebox(5.25,5.25){}}

\newtheorem{thm}{Theorem}
\newtheorem{pro}{Proposition}
\newtheorem{lem}{Lemma}

\newtheorem{corr}{Corollary}
\newtheorem{defin}{Definition}
\newtheorem{ob}{Observation}

\newenvironment{proposition}{\begin{pro} \nopagebreak}{{\hfill$\Box$} \end{pro}}
\newenvironment{thm-prf}{\begin{thm} \nopagebreak}{\end{thm}}
\newenvironment{pro-prf}{\begin{pro} \nopagebreak}{\end{pro}}
\newenvironment{lemma}{\begin{lem} \nopagebreak}{{\hfill$\Box$} \end{lem}}
\newenvironment{lem-prf}{\begin{lem} \nopagebreak}{\end{lem}}

\newenvironment{corr-prf}{\begin{corr} \nopagebreak}{\end{corr}}

\newcommand {\be}{\begin{equation}}
\newcommand {\ee}{\end{equation}}

\begin{document}

\title{Integrated Design of Unmanned Aerial Mobility Network: A Data-Driven Risk-Averse Approach}
\author[1]{Wenjuan Hou}
\author[2]{Tao Fang}
\author[2]{Zhi Pei}
\author[1]{Qiao-Chu He \thanks{Corresponding author: heqc@sustech.edu.cn} }
\affil[1]{School of Business, Southern University of Science and Technology}
\affil[2]{Department of Industrial Engineering, Zhejiang University of Technology}
\maketitle

\begin{abstract}

The real challenge in drone-logistics is to develop an economically-feasible Unmanned Aerial Mobility Network (UAMN). In this paper, we propose an integrated airport location (strategic decision) and routes planning (operational decision) optimization framework to minimize the total cost of the network, while guaranteeing flow constraints, capacity constraints, and electricity constraints. To facility expensive long-term infrastructure planning facing demand uncertainty, we develop a data-driven risk-averse two-stage stochastic optimization model based on the Wasserstein distance. We develop a reformulation technique which simplifies the worst-case expectation term in the original model, and obtain a fractable Min-Max solution procedure correspondingly. Using Lagrange multipliers, we successfully decompose decision variables and reduce the complexity of computation.

To provide managerial insights, we design specific numerical examples. For example, we find that the optimal network configuration is affected by the ``pooling effects" in channel capacities. A nice feature of our DRO framework is that the optimal network design is relatively robust under demand uncertainty. Interestingly, a candidate node without historical demand records can be chosen to locate an airport. We demonstrate the application of our model for a real medical resources transportation problem with our industry partner, collecting donated blood to a blood bank in Hangzhou, China.

\end{abstract}

\newpage
\section{Introduction}
Traditional logistics service based on land vehicles is facing emerging challenges such as increasing labor costs and traffic congestions, in contrast to growing logistics demand and rising customer expectation of delivery service quality. This is particularly the case for certain sectors in urban logistics, e.g., grocery or fresh products logistics, healthcare resources and products, and emergency/humanitarian logistics, in which context the goods are to be delivered within hours.
Drones logistics, as complement to land vehicles, can fill in the gaps in the aforementioned context and is already becoming a reality.
Although drone logistics cannot replace land vehicles entirely due to relatively higher cost and limited capacity, it certainly outperforms the latter in terms of efficiency (point-to-point delivery directly from origins to destinations) and accessibility (patching geographic holes such as islands and mountainous areas).
However, apart from certain showcase examples of drone delivery, it is not until recently that the technology (e.g., affordable AI-enabled navigation) and policy conditions \footnote {The Federal Aviation Administration (FAA) said on Friday, August 2 that it has approved the drone's first flight beyond the line-of-sight of the operator.} are prepared. Therefore, industry giants \footnote {At the end of March 2018, Amazon was patented that express drones can  recognize and respond to waving arms, pointing, flashing lights and voice, etc.} start prototype programs, but they quickly identify that the real challenge to implement economically-feasible drone-logistics at large scale, is to develop an unmanned aerial mobility network (UAMN): a data-and algorithm-driven service platform to coordinate the deployment of infrastructures (airports, charging stations, etc.) as well as the operations management of drones fleets (route planning, on-demand delivery scheduling,etc.).
Therefore, the industry needs motivate our research in designing such an UAMN, which is clearly relevant to the operations research/management community.

The existing literature on drone logistics pivots heavily on solving vehicle routing problems, with a few exceptions related to network design, e.g., \cite{boutilier2019response, Moon}. We depart from the on-demand scheduling/vehicle routing problems, and focus on route planning.
Thus, our model is most closely related to the hub-location model for air transportation network design \cite{shen2019reliable} and the reference therein.
The rationale is that aerial routes cannot be optimized on demand, since any route has to be reported to the regulating agency (such as FAA in the US) for prior approval. Once the routes are determined, the service platform can decide the frequency and volume en route but cannot adjust the fixed itineraries at will. In other words, the real problem of UAMN design is more analogous to building a highway network or urban bus network rather than vehicle routing problems.

In this paper, we address the following research question:
 how to design an UAMN under demand uncertainty while integrating macro-level strategic decisions (infrastructure planning decisions, e.g., airport location and capacity, as well as route planning) and micro-level operational decisions (transportation decisions, battery charging feasibility, etc.).
We choose a data-driven \emph{distributionally robust optimization} (DRO) framework to support the costly long-term infrastructure planning decisions due to demand uncertainty.

In this paper, we proposed a novel integrated airport location and routes planning model for UAMN, using a Wasserstein-distance based DRO framework. Our objective is to mininize the aggregate cost, with airport capacity constraints, flow-conservation constraints, battery capacity constraints, and satisfy all demands. We do not keep track of individual drones, but focus on transportation volumes (similar to bus frequency). We also abstract away from the dynamic charging process, and ensure that total electricity boosts can last the entire itineracy on a cyclic basis.

 We allow for mis-calibration of demand distribution by considering its ambiguity set based on the Wasserstain distance. Such Wasserstain distance-based distributionally robust optimization is flexible not only in fully utilizing historical demand data, but also risk-averse in the sense that inaccurate estimate for distribution leads to sub-optimal solution, e.g., \cite{zhao2018data}.
Such a decision framework is suitable for our problem context because costly long term infrastructure planning decisions are made under demand uncertainty with limited data.
We also reformulate the problem to maintain computational tractability. Numerical examples are proposed to generate managerial insights managerial insights concerning the strategic location decisions. Finally, we apply our model to solve a representative problem for our industry partner Xunyi\footnote{On 2019-10-17, our industry partner was issued a ``specific unmanned aerial vehicle commissioning letter'' and ``unmanned aerial vehicle logistics distribution business license''from the Civil Aviation Administration of China, the first approval of its kind in China (http://chinaplus.cri.cn/news/china/9/20191017/368089.html). In contrast, Amazon submitted FAA Approval for drone delivery on August 2019 and ``there is no set timeline for approval'' (https://www.aviationtoday.com/2019/08/09/following-wing-ups-amazon-seeks-approval-prime-air-drone-delivery/).} which is the first operating aerial logistics company approved by the Civil Aviation Administration of China.

Our specific contributions are:

1. We identify and conceptualize the key challenge in drone logistics as the UAMN design problem. Motivated by our industry partner, we focus on ``building aerial highway'' rather than solving vehicle routing problems.

2. We apply the Wasserstain distance-based distributionally robust optimization framework to facilitate data-driven decision making for both network design and operations management. We also develop computational technique in reformulating the complex problem into tractable models.

3. Through carefully designed numerical examples, we generate managerial insights towards the strategic infrastructure planning decisions.

4. Finally, in  collaboration with our industry partner, we apply our model to solve a real problem, delivering healthcare resources (blood and blood products) between blood collection station and the blood bank in densely populated urban areas. We summarize this application in a representative case study.

The rest of this paper is organized as follows. Section \ref{s-lit} reviews
relevant literature. Section~\ref{s-model} introduces our model setup. In
Section \ref{s-analysis}, we carry out the analysis. In Section \ref{s-num},
we describe our computational method, and provide numerical examples.
Section~\ref{s-con} concludes this paper with a discussion of future
research directions.
\section{Related Work}

\label{s-lit}

There is growing research interest in evaluating the economic efficiency in drone logistics, for example, in a hybrid truck-and-drone model, e.g.  \cite{agatz2018optimization, jeong2019truck, murray2015flying}. \cite{carlsson2017coordinated} demonstrated the improvement in efficiency this way, by combining a theoretical analysis in the Euclidean plane with real-time numerical simulations on a road network. \cite{poikonen2019branch} used a branch-and-bound approach to solve the TSP-D in which at each node the approximate lower bound is given by a dynamic program. \cite{poikonen2019mothership} are first to study the mothership and drone routing problem (MDRP) which is a more generalized TSP. However, the existing literature ignores the key challenge in drone-logistics, which is the UAMN design problem, so we focus our research on ``building aerial highway''.

The network design/facility location optimization have been studied extensively in the operations research community. \cite{shen2019reliable} study a reliable hub location model, by exploiting the structural properties of the problem, they introduce a tractable mixed-integer linear program reformulation and develop a constraint generation method to accelerate the solution procedure.
\cite{an2014reliable} develop a two-stage RO reliable p-median facility location model to minimize the weighted cost in normal and disruption scenarios. \cite{chowdhury2017drones, cui2010reliable} develop a continuum approximation (CA) model aiming to minimize the overall expected total cost in normal and failure/disruption scenarios. Based on \cite{an2014reliable, chowdhury2017drones, cui2010reliable}, \cite{lim2013facility} develop a stylized continuous model, and investigate the impact of misestimating the disruption probability.
\cite{lu2015reliable} present a model that allows disruptions to be correlated with an uncertain joint distribution. And a distributionally robust optimization is applied to minimize the expected cost under the worst-case distribution with given marginal disruption probabilities. \cite{mak2013infrastructure} aid to locate battery-swapping infrastructure and choose its capacity along an existing network of freeways with by using limited information of demand, e.g., mean, variance, and develop a chance-constraint robust optimization model. In this paper, we solve the integrated problem of both facility location and network design, fully utilizing the historical data, and develop a Wasserstain-distance robust optimization model.

In addition, UAMN design problems can also get some inspiration from the traditional airline operation management.
\cite{dunbar2012robust} introduces a new approach to accurately calculate and minimize the cost of propagated delay in a framework that integrates aircraft routing and crew pairing.
 \cite{DUNBAR201468} extends the approach of \cite{dunbar2012robust} by proposing two new algorithms that achieve further improvements in delay propagation reduction via the incorporation of stochastic delay information. \cite{WEI201413} introduces the flight routes addition/deletion problem and compares three different methods to analyze and optimize the algebraic connectivity of the air transportation network. \cite{sun2008multicommodity} aims to manage the increasing complexity of traffic flow in the airspace and present a traffic flow model called the Large-Capacity Cell Transmission Model in which the integer program is relaxed to a linear program for computational efficiency. \cite{WEI20131} inherits this problem, rewrite it in a standard linear programming and analyze the total unimodular property of the constraint matrix. The authors prove that a simplex related method guarantees the solutions to be optimal.
 \cite{doi:10.2514/1.I010710} focuses on the formulation of fixed finial time multiphase optimal control problem with energy consumption as the performance index for a multirotor eVTOL aircraft. \cite{inproceedings} proposes mathematical models for commercial transport service providers to decide which type of schedule to offer, how to dispatch the fleet and schedule operations, based on simulated market demand, such that profits is maximized.
However, this stream of research usually step from a small segmentation of aerial mobility network and develop a heuristic approach. To generate systematic conclusions and managerial insight, our basic model is somewhat stylized and can be solved via an exact approach.

Stochastic programming can effectively describe many decision-making problems in uncertain environments. We finally review the solution methodology for stochastic programming in the network design and facility location literature. We do not post distributional assumptions on demand uncertatinty (typically required for stochastic programming with chance constraints). Note that the solution framework we used in this paper has been also applied to the unit commitment problems \cite{wang2011chance, wu2014chance}. The classic solution to such stochastic program is by scenario sampling \cite{Calafiore2005}, or robust optimization \cite{mak2013infrastructure}. \cite{JIANG2014751} proposed a two-stage robust optimization model to address the network constrained unit commitment problem under uncertainty.
\cite{mete2010stochastic} studies the storage and distribution problem of medical supplies by implementing stochastic optimization approach.
In this paper, our solution approach is most related to \cite{zhao2018data} which studied a data-driven risk-averse stochastic optimization approach with Wasserstein distance for the general distribution case. They reformulated the risk-averse two-stage stochastic optimization problem to a traditional two-stage robust optimization problem via using Wasserstein distance.
We apply this framework to our UAMN design problem, and develop an integrated data-driven risk-verse model.

\section{Model}

\label{s-model}

\textbf{Demand uncertainty.}
We use $b_k(\omega )$ to denote the demand volume for the kth O-D pair.  Considering that there is uncertainty in the demand, $b_{k}(\omega )$ is a random variable, wherein each $\omega$ corresponds to a random draw from some sample space $\Omega $. We can omit $\omega$ in $b_{k}(\omega )$ and denote demand volume by $b_{k}$.
In fact, the actual demand distribution for $b_k$  is unknown at the stage of infrastructure planning, e.g., airport location, route and channel design, etc. An inaccurate estimate of the demand distribution may lead to suboptimal infrastructure investment decisions, and thus, it is desirable to adopt a distributionally robust optimization framework without picking a fixed distribution ex ante.
At this stage, we have no restriction for the demand distribution other than that it is bounded with an upper bound (denoted by $W_k^+$) and a lower bound (denoted by $W_k^-$).
Meanwhile, to utilize the historical demand data, we adopt a data-driven approach wherein the true distribution should be ``close'' to the reference distribution empirically estimated from data, in the sense of ``Wasserstein distance metric'' which has been utilized in \cite{zhao2018data}, \cite{chen}. We will specify the detailed setup in the analysis.

\textbf{Unmanned aerial mobility network.}
We use a set $V$ to denote all nodes in the UAMN, which includes supply nodes (``origin'', i.e., warehouses and distribution centers), demand nodes (``destination'', i.e., customer) and transfer nodes (i.e., transfer airport).
Goods are delivered from origins to destinations though the UAMN, enabled by a fleet of orchestrated drones via optimization algorithms.
We represent a delivery demand $k$ by its origin $O_k$ and destination $D_k$, and use a set  $K$ to denote all O-D pairs.
A pair of nodes in the UAMN constitute a route, represented by the arc set $A$ in the network.
A ``route'' connecting two nodes consists of multiple ``channels'', i.e., aerial corridors that confines the trajectories of drones, similar to a highway consisting of multiple lanes. In addition, there could be multiple types of channels (fast vs. slow) and we denote the channel type by a set $T$.
A node with a non-zero demand indicates an airport, at which drones can pick up or deliver goods, as well as lay over for charging.
A node is referred to as a ``transfer airport'', if drones can charge but not pick up or deliver at the node.
We use a binary decision variable $z_i$ to decide whether an airport should be in potential location $i$.
In this paper, we focus on the UAMN design and do not keep track of the moving path of individual drones.

\textbf{Flow constraints.}
For a particular O-D pair (denoted by index $k$), we use set $V-O_k-D_k$ to denote its transfer nodes. We represent the demand volume of $k$ by $b_k$.
We use $x_{ji}^{k}$ to denote transportation volume en route $(i,j)$ for kth O-D pair.
For a particular O-D pair (denoted by index $k$), at a pick-up location (origin node), the net transportation volume out of the node is $b_k$ (flow-out volume subtracted by the flow-in volume):
\begin{equation*}
\sum_{j:(i,j)\in A}x_{ij}^{k}(\omega)-\sum_{j:(i,j)\in A}x_{ji}^{k}(\omega)=b_{k}(\omega),\forall
k\in K,\forall i\in O_{k}.
\end{equation*}%
At a pick-up location (demand node), the net transportation volume out of the node is $-b_k$ (flow-out volume subtracted by the flow-in volume):
\begin{equation*}
\sum_{j:(i,j)\in A}x_{ij}^{k}(\omega)-\sum_{j:(i,j)\in A}x_{ji}^{k}(\omega)=-b_{k}(\omega),\forall
k\in K,\forall i\in D_{k}.
\end{equation*}%
At a pick-up location (transfer node), the net transportation volume out of the node is $0$ (flow-out volume subtracted by the flow-in volume):
\begin{equation*}
\sum_{j:(i,j)\in A}x_{ij}^{k}(\omega)-\sum_{j:(i,j)\in A}x_{ji}^{k}(\omega)=0,\forall k\in
K,\forall i\in V-O_{k}-D_{k}.
\end{equation*}%

\textbf{Capacity constraints.}
We first make strategic decisions by locating airports. An airport is capacitated due to (1) limited ``port'' for landing and take-off as well as (2) limited parking spaces, and the capacity for the airport at location $i$ is $w_i$.
As we have described, a route may consist of multiple channels of type $t$, we assume that the capacity of a type t channel en route $(i,i)$ is $u_{ij}^t$.
 We use an integer decision variable $y_{ij}^t$ to decide how many channels of type $t$ should be invested in a potential route $(i,j)$.

For a route denoted by an arc $(i,j)$ in the UAMN,
the total net transportation volume through this route cannot exceed total capacity volume.  Hence, we set this constraint to limit transportation volume. The left side represents the sum of transportation volume over all O-D pairs, while the right hand represents the sum of capacity volume over all channel types.
\begin{equation*}
\sum_{k\in K}x_{ij}^{k}(\omega)\leq \sum_{t\in T}u_{ij}^{t}y_{ij}^{t},\forall
(i,j)\in A.
\end{equation*}%

For a location node denoted by $i$, if $z_i=1$ (an airport should be), its capacity volume is $w_{i}$, if $z_i=0$ (an airport should not be), its capacity is unlimited. We integrate these two possible cases together and set the constraint as follows, wherein $M$ is a $large$ number:
\begin{equation*}
\sum_{k\in K}\sum_{j:(i,j)\in A}x_{ij}^{k}(\omega)+\sum_{k\in K}\sum_{j:(i,j)\in
A}x_{ji}^{k}(\omega)\leq w_{i}z_{i}+M(1-z_{i}),\forall i\in V.
\end{equation*}

Notice that the capacity constraints both at the airport and en route between nodes relate the location decisions $\{z\}$'s, network design decisions $\{y\}$'s, and the transportation decisions $\{x\}$'s. Therefore, the strategic (infrastructure) decisions are intertwined with the operation decisions (transportation or scheduling). This is a salient feature of UAMN design in terms of the integration of different decision hierarchies.

\textbf{Battery capacity constraints.}
We use $l_{ij}$ to denote the consumption of electricity quantity en route $(i,j)$, while $L$ denotes a lumpsum boost in battery level upon a charging completion.
For a particular O-D pair (denoted by index $k$), the total battery consumption is less than the total battery charged. We set the electricity constraints as follows, wherein the first term is the sum of battery consumption over all routes, and the second term is the sum of electricity charged over all routes. The amount of electricity charged is only included in the summation when the route is chosen ($x\neq0$) and the airport (charging station) is located ($z=1$). This constraint means that there is enough total charge to reach the charge-discharge balance in one operating period.
\begin{equation*}
\sum_{(i,j)\in A}l_{ij}x_{ij}^{k}(\omega)\leq b_k\sum_{(i,j)\in A}z_{i}L,\forall
k\in K,
\end{equation*}

We also abstract away from the dynamic charging process to focus on the charging location: We assume that drones can only get charged at airports (including transfer ones), and we ensure that total electricity boosts (can charge multiple times) can last the entire itineracy, on a cyclic basis.
As we focus on strategic level planning, we do not keep track of the battery level for individual drones (also computationally intractable).

\textbf{Objective.}
The aggregate costs consist of transportation cost, aerial channel cost, airport infrastructure cost, and airport capacity cost. We use $Ct_{ij}^k$ to denote the cost to transport unit volume cargo for the kth O-D pair en route $(i,j)$, $Cd_{ij}^t$ to denote the cost to set up a channel of type $t$ en route $(i,j)$, $Cf_j$ to denote the infrastructure cost to locate an airport for the jth location, and $Cs_j$ to denote the cost to set unit volume capacity for the jth location.
We sum up all these costs over all routes, all locations, and all types. Our objective function takes the form as follows:
\begin{equation*}
\min \sum_{k\in K}\sum_{(i,j)\in A}Ct_{ij}^{k}x_{ij}^{k}(\omega)+\sum_{t\in
T}\sum_{(i,j)\in A}Cd_{ij}^{t}y_{ij}^{t}+\sum_{j\in V}Cf_{j}z_{j}+\sum_{j\in
V}Cs_{j}w_{j}z_{j}.
\end{equation*}%
Our model is systematic by integrating all kinds of costs.

To summarize, we solve the following mixed integer program:
\begin{itemize}
  \item Objective:
\begin{equation*}
\min \sum_{k\in K}\sum_{(i,j)\in A}Ct_{ij}^{k}x_{ij}^{k}(\omega)+\sum_{t\in
T}\sum_{(i,j)\in A}Cd_{ij}^{t}y_{ij}^{t}+\sum_{j\in V}Cf_{j}z_{j}+\sum_{j\in
V}Cs_{j}w_{j}z_{j}.
\end{equation*}%
  \item Flow constraints:
\begin{equation}
\sum_{j:(i,j)\in A}x_{ij}^{k}(\omega)-\sum_{j:(i,j)\in A}x_{ji}^{k}(\omega)=b_{k}(\omega),\forall
k\in K,\forall i\in O_{k}.
\end{equation}%

\begin{equation}
\sum_{j:(i,j)\in A}x_{ij}^{k}(\omega)-\sum_{j:(i,j)\in A}x_{ji}^{k}(\omega)=-b_{k}(\omega),\forall
k\in K,\forall i\in D_{k}.
\end{equation}%

\begin{equation}
\sum_{j:(i,j)\in A}x_{ij}^{k}(\omega)-\sum_{j:(i,j)\in A}x_{ji}^{k}(\omega)=0,\forall k\in
K,\forall i\in V-O_{k}-D_{k}.
\end{equation}%
  \item Capacity constraints:
\begin{equation}
\sum_{k\in K}x_{ij}^{k}(\omega)\leq \sum_{t\in T}u_{ij}^{t}y_{ij}^{t},\forall
(i,j)\in A.
\end{equation}%

\begin{equation}
\sum_{k\in K}\sum_{j:(i,j)\in A}x_{ij}^{k}(\omega)+\sum_{k\in K}\sum_{j:(i,j)\in
A}x_{ji}^{k}(\omega)\leq w_{i}z_{i}+M(1-z_{i}),\forall i\in V.
\end{equation}
  \item Battery capacity constraints:
\begin{equation}
\sum_{(i,j)\in A}l_{ij}x_{ij}^{k}(\omega)\leq b_k\sum_{(i,j)\in A}z_{i}L,\forall
k\in K,
\end{equation}
\end{itemize}

\bigskip
We summarize the nomenclature in this paper as follows:
\begin{table}[H]
\begin{center}
\begin{tabular}{l p{190pt} l p{190pt}}
  \hline
  $A$ & The set of all arcs in the network. &
  $V$ & The set of all nodes in the network. \\
  $K$ & The set of all O-D pairs in the network. &
  $T$ & The set of all types of channels (e.g., fast, slow, etc.) in the arcs. \\
  $O_{k}$ & The origin node for an O-D pair.&
  $D_{k}$ & The destination node for an O-D pair.\\
  $b_{k}$ & The transportation demand for an O-D pair. &
  $x_{ij}^{k}$ & Continuous variables for transportation volume. \\
  $y_{ij}^{t}$ & Integer variables for number of channel of type $t$ at route $(i,j)$.&
  $u_{ij}^{t}$ & Capacities of channel of type $t$ at route $(i,j)$. \\
  $z_{i}$ & Binary decision variables for an airport at location $i$. &
  $w_{i}$ & Capacities for an airport at location $i$. \\
  $L$ & A lumpsum boost in battery level upon a charging completion.&
  $l_{ij}$ & Energy/electricity consumed for route $(i,j)$.\\
  $Ct_{ij}^k$ & The cost to transport unit volume cargo for the kth O-D pair en route $(i,j)$.  &
  $Cd_{ij}^t$ & The cost to set up a channel of type $t$ en route $(i,j)$.\\
  $Cf_j$ & The infrastructure cost to locate an airport for the jth location. &
  $Cs_j$ & The cost to set unit volume capacity for the jth location.\\
  \hline
\end{tabular}
\end{center}
\end{table}
\bigskip
\section{Analysis and Computation Technique}

\label{s-analysis}

\subsection{Reformulations}

We set up a two-stage model in last section, in which the first-stage makes decisions for locations and channels and the second-stage makes decisions for transformation volume. Its second-stage problem is as follows:
\begin{align*}
\mathcal{Q}(\eta,b)=\min_{x}&\quad \sum_{k\in K}\sum_{(i,j)\in A}Ct_{ij}^{k}x_{ij}^{k}\\
\text{s.t.}\quad &(1),(2),(3),(4),(5),(6)
\end{align*}%
where $\eta=(z,y)$. We use $\eta$ to vectorize decision variables of $z$ and $y$.

We can rewrite the $\mathcal{Q}(\eta,b)$ as follows:
   \begin{align*}
        \min_{X}\quad &C^{T}X\\
        \text{s.t.}\quad &AX=B\\
        &DX\leq E,
   \end{align*}
Where $A$, $B$, $C$, $D$, $E$ are matric, and we give the detail definition of them in appendix.

Using traditional two-stage stochastic optimization framework in \cite{A.Shapiro}, we can formulate above two-stage problem as follows:
\begin{equation*}
(SP)\min_{\eta}E_{P}[\mathcal{Q}(\eta,b)]+\sum_{t\in
T}\sum_{(i,j)\in A}Cd_{ij}^{t}y_{ij}^{t}+\sum_{j\in V}Cf_{j}z_{j}+\sum_{j\in
V}Cs_{j}w_{j}z_{j}
\end{equation*}%
where $P$ is the actual distribution of demand. However, $P$ is usually unknown and the estimate of it is different and inaccurate, i.e., historical demand data are not enough, future demand may change. In addition, this model incorporates no risk averseness.

Considering that the cost to construct an new UAMN is expensive, the UAMN needs to have a long-term utility under demand uncertainty. We choose the worst-case distribution in the set defined by Wasserstein distance metric and consider a data-driven risk-averse stochastic optimization formulation as follows:
\begin{equation*}
(DD-SP)\min_{\eta} \max_{\widehat{P}\in \mathcal{D}}E_{\widehat{P}}[\mathcal{Q}(\eta,b)]+\sum_{t\in
T}\sum_{(i,j)\in A}Cd_{ij}^{t}y_{ij}^{t}+\sum_{j\in V}Cf_{j}z_{j}+\sum_{j\in
V}Cs_{j}w_{j}z_{j}
\end{equation*}%
\begin{equation*}
\mathcal{D}=\{\widehat{P}:d_{M}(P_{0},\widehat{P})\leq\theta\}
\end{equation*}%
where $P_{0}$ is the reference empirical distribution determined through historical data. $d_{M}$ is the Wasserstein distance defined on two distributions. $\mathcal{D}$ is a confidence set of the true distribution by confining $d_{M}$ less than or equal to $\theta$ which decides the size of $\mathcal{D}$. On the one hand, we fully utilize the historical demand data because $\mathcal{D}$ is a set near the reference distribution. On the other hand, we use a set to estimate actual distribution rather than a concrete distribution. So our model is both data-driven and risk-averse.

The aforementioned decision framework is known to have verified convergence property \cite{zhao2018data}. However, directly solving this two-stage risk-averse stochastic program is computationally challenging. Next, we start by reformulating the problem.

\begin{proposition}
The problem $DD-SP$ under $\mathcal{D}$ is equivalent to the following two-stage robust optimization problem:
\begin{equation*}
(RDD-SP)\min_{\eta,\beta\geq0}\frac{1}{N}\sum_{i=1}^{N}\max_{b}[{\mathcal{Q}(\eta,b)-\beta\rho^{i}(b)}]+\theta\beta+\sum_{t\in
T}\sum_{(i,j)\in A}Cd_{ij}^{t}y_{ij}^{t}+\sum_{j\in V}Cf_{j}z_{j}+\sum_{j\in
V}Cs_{j}w_{j}z_{j}
\end{equation*}%
where $\rho^{i}(b)=\rho(b,b^i)$.
\end{proposition}

We reformulate the original problem $(DD-SP)$ to a two-stage robust optimization problem in Proposition 1. By finding the worst distribution of $b$ in $\mathcal{D}$, we eliminate the confidence set and attain the expected value. The worst distribution of $b$ is also discrete because Wasserstein distance is used to capture the optimal transport between discrete distributions. The term $\beta\rho^{i}(b)$ can be understood as a penalty about $b$. The value of $\beta$ decides the confidence level of $\mathcal{D}$. The other terms remain the same as that in the original problem $(DD-SP)$. However, the subproblem is still not linear because $\rho$ is nonlinear, which has the absolute value. We will deal with $\rho$ in the next proposition.

\begin{proposition}
The problem $RDD-SP$ is equivalent to the following Min-Max problem:
\begin{align}
(Master)\min_{\eta,\beta\geq0}&\frac{1}{N}\sum_{j=1}^{N}\varpi^j(\eta,\beta)+\theta\beta+\sum_{t\in
T}\sum_{(i,j)\in A}Cd_{ij}^{t}y_{ij}^{t}+\sum_{j\in V}Cf_{j}z_{j}+\sum_{j\in
V}Cs_{j}w_{j}z_{j}\notag
\end{align}%
where $\eta=(\omega,z,y)$, and $\varpi^j(\eta,\beta)$ is a (SUB) problem as follows:
\begin{align*}
(Sub)\varpi^j(\eta,\beta)
=&\max_{\delta^{+},\delta^{-},\mu,\lambda\geq0}\quad [-\sum_{i=1}^{K}(\mu_{i}-\mu_{K+i})((W_i^{+}-b_i^j)\delta_i^{+}+(W_i^{-}-b_i^j)\delta_i^{-}+b_i^j)\\
&-\beta\sum_{i=1}^K((W_i^{+}-b_i^j)\delta_i^{+}-(W_i^{-}-b_i^j)\delta_i^{-})-\lambda^TE]\\
\text{s.t.}\quad  &C^T+\mu^TA+\lambda^TD\geq0\\
&\delta_i^{+},\delta_i^{-}\in\{0,1\}\\
&\delta_i^{+}+\delta_i^{-}\leq1.
\end{align*}
where $W_i^{+}$ is the upper bound of $b_i$, and $W_i^{-}$ is the lower bound of $b_i$.
\end{proposition}

The salient difference between Proposition 1 and Proposition 2 is that we reformulate the subproblem to a linear problem for a given first stage decision. We attain this subproblem by using the dual form of $\mathcal{Q}(\eta,b)$, and linearizing $\rho^{i}(b)$. This Min-Max problem is a traditional two-stage robust problem£¬ in which the first stage is the master problem and the second stage is the subproblem.

Because we linearize the subproblem in proposition 2, this Min-Max problem is solvable. For any $y$, $z$ and $\beta$, we can calculate the corresponding optimal value of the subproblem. Hence, we can calculate the optimal value of the master problem and optimal solutions to the decision variables $y$ and $z$.
Because $\beta$ is the unique optimal value for every $\theta$, it suffices to obtain the value of $b$ by solving the subproblem for given $\beta$. The confidence interval (radius of the Wasserstein ball) can be adjusted by choosing $\beta$.

\subsection{Solution procedures}

In the previous subsection, we have reformulated our model as a Min-Max problem. Because the decision variables $z$ and $y$ occur both in the master problem and subproblem, they have to be solved recursively. For given choices of $y$ and $z$, the subproblem will be a tractable linear program. Then we solve the entire integer optimization problem with respect to $y$ and $z$. We describe the solution procedures in this section with a sketch of the derivations. The readers are referred to Appendix D for detailed derivations.

\textbf{Step1.} We notice that the matric $E$ in the subproblem is composed of decision variables $y$ and $z$, so we decompose $E$. Then we collect terms with respect to $y$ and $z$.

\textbf{Step2.} We notice that these decision variables $y$ and $z$  rely on $\mu$ and $\lambda$, respectively. However, decision variables $y$ and $z$ are interdependent because constraints are comprised of $\mu$ and $\lambda$. By using multiplier $\gamma\geq0$, we can eliminate the constraint. We can rewrite the master problem as follows:
\begin{align}
\psi=&\max_{\delta^{+},\delta^-,\mu,\lambda\geq0}\min_{\eta,\beta\geq0}\quad -\frac{1}{N}\sum_{j=1}^{N}
[\sum_{i=1}^{K}(\mu_{i}^j-\mu_{K+i}^j)((W_i^{+}-b_i^j)\delta_i^{j+}+(W_i^{-}-b_i^j)\delta_i^{j-}+b_i^j)]\\
+&[-\frac{1}{N}\sum_{j=1}^{N}
\sum_{i=1}^K((W_i^{+}-b_i^j)\delta_i^{j+}-(W_i^{-}-b_i^j)\delta_i^{j-})+\theta]\beta\notag\\
+&\sum_{j\in V}[\{-(\overline{\lambda}_{1})_j(w_j-M)+Cf_{j}+Cs_{j}w_{j}\}z_j-(\overline{\lambda}_{1})_jM]\notag\\
+&\sum_{t\in T}\sum_{(i,j)\in A}\{Cd_{ij}^{t}-(\overline{\lambda}_{2})_{ij}u_{ij}^t\}y_{ij}^{t}\notag\\
+&\sum_{j=1}^{N}\{C^T+(\mu^j)^TA+(\lambda^j)^TD\}\gamma^j\notag\\
\text{s.t.}\quad
&\delta_i^{+},\delta_i^{-}\in\{0,1\}\notag\\
&\delta_i^{+}+\delta_i^{-}\leq1.\notag
\end{align}
where $\lambda_1$ represents the multiplier of constraints (5), $\lambda_2$ represents the multiplier of constraints (4), and $\lambda_3$ represents the multiplier of constraints (6).  So the multiplier $\lambda$ can be separated by $\lambda=({\lambda_2}^T,{\lambda_1}^T,{\lambda_3}^T)^T$. $\lambda_1^j$ is corresponding to $j$th sample of $\lambda_1$, $\lambda_2^j$ is corresponding to $j$th sample of $\lambda_2$, and $\lambda_3^j$ is corresponding to $j$th sample of $\lambda_3$.
 $\overline{\lambda}_1$ is the mean of $\lambda_1^j$  , and $\overline{\lambda}_2$ is the mean of $\lambda_2^j$. So the multiplier $\lambda^j$ can be separated by $\lambda^j=({\lambda_2^j}^T,{\lambda_1^j}^T,{\lambda_3^j}^T)^T$.
$(\overline{\lambda}_{1})_j$ is the component in $\overline{\lambda}_{1}$ corresponding to $z_j$. $(\overline{\lambda}_{2})_{ij}$ is the component in $\overline{\lambda}_{2}$ corresponding to $y_{ij}$.

\textbf{Step3.} We have eliminated the constraint by the above operation. We can notice that $z$ and $y$ are related to the different part of $\lambda$. If we can successfully separate $\lambda$, then the decision of $z$ and the decision of $y$ will be independent. By separating $\sum_{j=1}^{N}\{C^T+(\mu^j)^TA+(\lambda^j)^TD\}\gamma^j$, we rewrite the problem as follows:
\begin{align*}
\psi
=&\max_{\delta^{+},\delta^-,\mu}\min_{\beta\geq0}\{-\frac{1}{N}\sum_{j=1}^{N}
[\sum_{i=1}^{K}(\mu_{i}^j-\mu_{K+i}^j)((W_i^{+}-b_i^j)\delta_i^{j+}+(W_i^{-}-b_i^j)\delta_i^{j-}+b_i^j)-(\mu^j)^TA\gamma^j]\tag{8}\\
+&[-\frac{1}{N}\sum_{j=1}^{N}
\sum_{i=1}^K((W_i^{+}-b_i^j)\delta_i^{j+}-(W_i^{-}-b_i^j)\delta_i^{j-})+\theta]\beta\}\\
+&\max_{\lambda_1,\lambda_3\geq0}\min_{z}\{
\sum_{j\in V}[\{-(\overline{\lambda}_{1})_j(w_j-M)+Cf_{j}+Cs_{j}w_{j}\}z_j-(\overline{\lambda}_{1})_jM]+\sum_{j=1}^{N}(\lambda_1^j)^TD_1\gamma^j
+\sum_{j=1}^{N}(\lambda_3^j)^TD_3\gamma^j\}\\
+&\sum_{(i,j)\in A}\max_{\lambda_{2ij}\geq0}\min_{y_{ij}^t}\{
[\sum_{t\in T}\{Cd_{ij}^{t}-(\overline{\lambda}_{2})_{ij}u_{ij}^t\}y_{ij}^{t}+\sum_{h=1}^{N}(\lambda_{2}^h)_{ij}^T(D_{2})_{(ij)}\gamma^h]\}\notag\\
+&\sum_{j=1}^{N}C^T\gamma^j\\
\text{s.t.}\quad
&\delta_i^{+},\delta_i^{-}\in\{0,1\}\\
&\delta_i^{+}+\delta_i^{-}\leq1.
\end{align*}
where $D_1$ represents the coefficients in constraints (5), $D_2$ represents the coefficients in constraints (4), and $D_3$ represents the coefficients in constraints (6).  So the multiplier $D$ can be separated by $D=({D_2}^T,{D_1}^T,{D_3}^T)^T$. $(D_{2})_{(ij)}$ is a row of $D_2$ corresponding to $y_{ij}$.

\textbf{Step4.} We have decomposed $\lambda$ and $D$ according to $z$ and $y$. Given the proper multiplier value of $\gamma$, solving decision variables of $z$ and $y$ are independent. Furthermore, we can separate $y$ to $y_{ij}$, because every $y_{ij}$ is only related to $(\lambda_2)_{ij}$. We change the arrangement of sum and Max-Min, and attain the formula of (8). We can notice that the original Min-Max problemd which contains all decision variables has been divided into several independent Min-Max problem which contains a smaller number of variables. We notice that $y$ can be divided into $y_{ij}$ because every $y_{ij}$ is only related to $(\lambda_2)_{ij}$, but $z$ can not be divided into $z_j$ because $z_j$ is related to $\lambda_1$.  Given the value of $\gamma$, every decision variable $y_{ij}$ can be solved independently. In addition, given the value of $\gamma$, decision variables $z$ can be solved.

\textbf{Step5.} We must remember that the multiplier $\gamma$ represents the condition $C^T+(\mu^j)^TA+(\lambda^j)^TD\geq0$.
So we just need to find a proper  $\gamma$ make this constraint hold, then the remained computation is decreased dramatically, allowing us to conduct our numerical study later.

 It is clear that any node with nonzero transportation volume should be an airport location. We choose to check this constraint at the end of the algorithm to reduce the computation complexity.
 The formula (8) takes advantage of tractability comparing to the Min-Max problem which contains all decision variables. For the original Min-Max problem, the complexity increases exponentially in the size of $z$ and $y$. Though $z$ can not be divided into $z_j$, the amount of computation for (8) is dramatically decreased comparing to the Min-Max problem in last section because the size of $y$ is much larger than the size of $z$.

\section{Numerical Examples and Case Study}

\label{s-num}
\subsection{Numerical examples}
We propose numerical examples to generate managerial insights. In particular, we design the example to isolate the effect of particular factors on the network configurations.

To streamline our numerical examples, we fix some basic model parameters throughout this section unless specified otherwise.
We need to give value to some important parameters, the node set $V$, the confidence set parameter $\theta$, the confidence level $\beta=1000$, the number of channel type $t$, and the range of channel decision variable $y$.
we choose $V=5$, $\theta=100$, $\beta=1000$, $N=2$, $t=1$, $y$ is binary. The rest of the parameters are also carefully chosen within a reasonable range of values.
It suffices to choose $N=2$ to deliver the main messages in the following numerical examples. Demands are simulated by censored Gaussian distribution, with the upper bound and lower bound chosen at triple the standard deviations above and below the mean, respectively. Similarly, we carefully assign reasonable values to the location of nodes (and thus the pairwise distance between them), the length of routes, as well as cost coefficients specified in Table 1.

\begin{proposition}
(Observation 1) The marginal decision of setting up an additional airport at a particular candidate location depends crucially on both the infrastructure cost at and transportation cost en route to the location.
\end{proposition}
In this numerical example, we fix all parameters and vary solely the parameters listed in Table 1, one at a time.
Not surprisingly, when we either increase the infrastructure cost at or transportation cost en route to node 3, the optimality of airport location is to avoid node 3 and choose node 4.
This observation is intuitive. A direct implication is to decide whether to set up airport in densely populated urban areas (low transportation costs due to proximity towards demands). Nevertheless, it reflects a trade-off between infrastructure costs and transportation costs in locating the airport, which also justifies the integration of two dimensions (location vs. transportation) in the first place.
\begin{table}[H]
\caption{Transfer airport location}
\begin{center}
\begin{tabular}{l l l}
  \hline
  Optimal location for transit airport & Node 3 & Node 4 \\
  \hline
  Airport infrastructure cost & $Cf_4$=15000 & $Cf_4$=8000 \\
  Airport capacity cost & $Cs_4$=2 & $Cs_4$=1 \\
  Channel cost & $Cd_{i4}$=2000 & $Cd_{i4}$=1200 \\
  Transportation cost & $Ct_{i4}$=12 & $Ct_{i4}$=10\\
  \hline
\end{tabular}
\end{center}
\end{table}

\begin{figure}[H]
 \centering
 \includegraphics[width=0.80\textwidth]{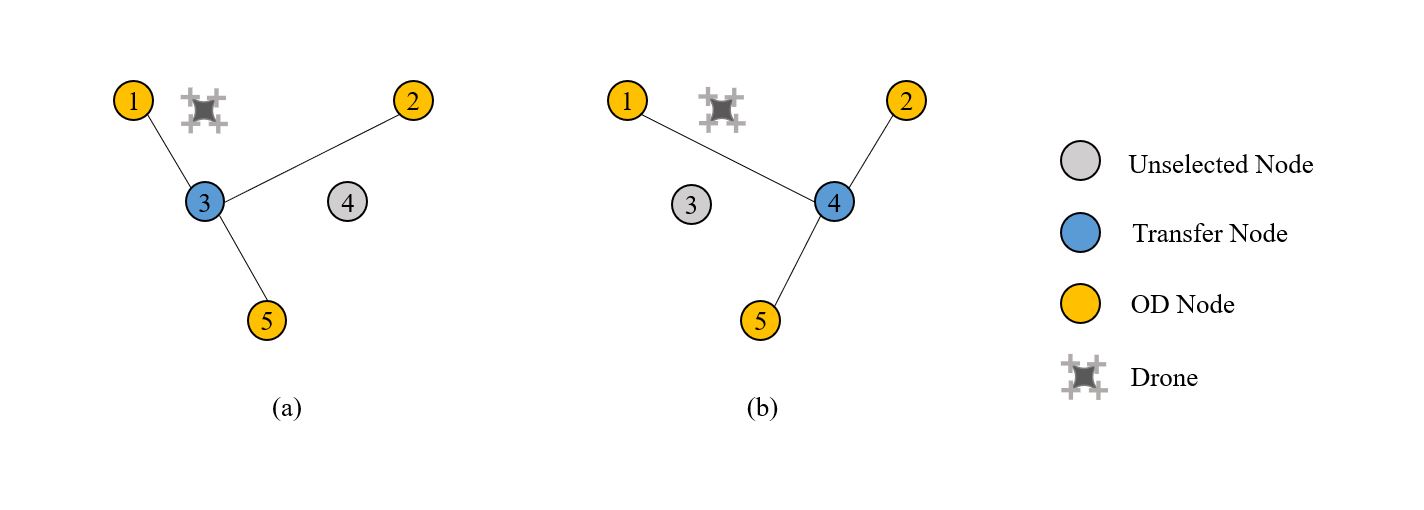}
 \caption{Transfer airport location}
\end{figure}
\begin{proposition} (Observation 2)
The optimal network configuration depends jointly on infrastructural and transportation cost coefficients,  the battery capacity constraints, as well as the pooling effects in channel capacities.
\end{proposition}

In this numerical example, we fix all parameters and vary solely the parameters listed in Table 2, one at a time.
When we either increase the channel cost at or transportation cost en routes between nodes 1,2,5, the optimality of topology is to avoid triangle configuration and chooses hub-and-spoken configuration. When we decrease battery capacity, it chooses hub-and-spoken configuration. In particular, when battery capacity is small, the triangle configuration can never be the optimal solution even if its cost is lower, because  drones need to satisfy the charge-discharge balance constraints.
Finally, we briefly explain the pooling effects in channel capacities. In general, we prefer a single high capacity channel rather than multiple low capacity channels due to economies of scale. The centralized channel smooths out the demand uncertainty and increases the utilization rate of the channel in general.

\begin{table}[H]
\caption{Topology}
\begin{center}
\begin{tabular}{l l l}
  \hline
  Optimal network design & Triangle configuration & Hub-and-spoken configuration \\
  \hline
  Channel cost  & $Cd_{12}$, $Cd_{15}$, $Cd_{25}$=1200 & $Cd_{12}$, $Cd_{15}$, $Cd_{25}$=7000\\
  Transportation cost  & $Ct_{12}$, $Ct_{15}$, $Ct_{25}$=10 & $Ct_{12}$, $Ct_{15}$, $Ct_{25}$=30\\
  Battery capacity  & $L=20$ & $L=10$\\
  \hline
\end{tabular}
\end{center}
\end{table}

\begin{figure}[H]
 \centering
 \includegraphics[width=0.80\textwidth]{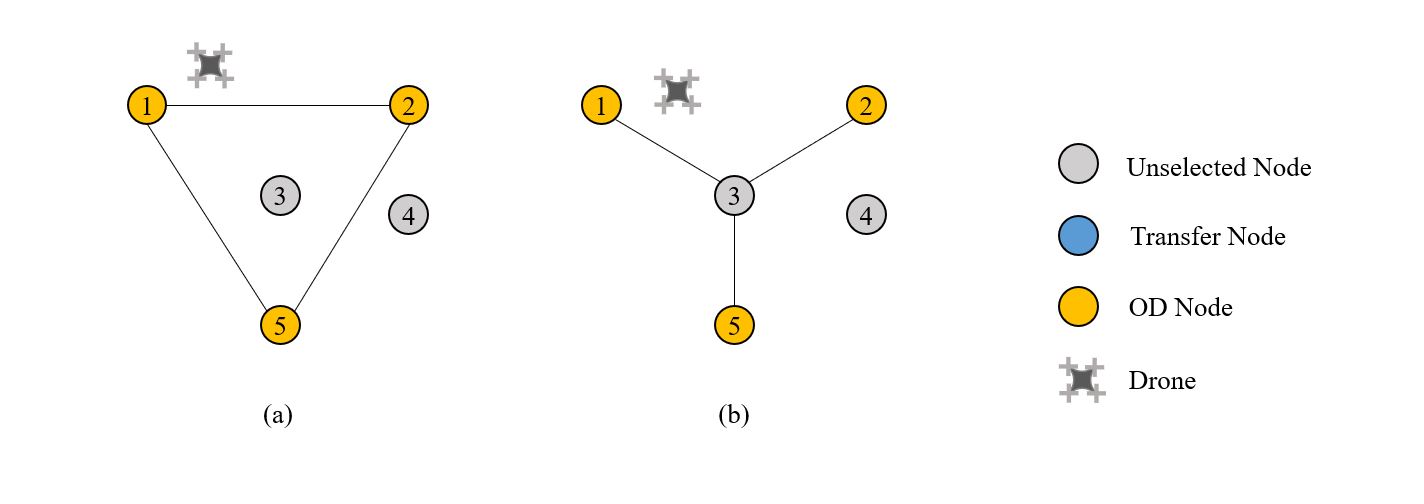}
 \caption{Topology}
\end{figure}
\begin{proposition} (Observation 3)
The optimal network design is relatively robust under demand uncertainty.
\end{proposition}
Notice that to switch from node 3 to node 4, we need to either increase the mean of demand by six times or the standard deviation by ten times. In other words, the network configuration is less sensitive towards the demand uncertainty, compared with other cost coefficients.
This observation echoes our choices of distributionally robust optimization framework to avoid suboptimal long-term infrastructure and investment under a conservative and adversarial demand forecast.

\begin{table}[H]
\caption{Effect of mean and variance}
\begin{center}
\begin{tabular}{c c c}
  \hline
  Optimal location for transit airport & Node 3 & Node 4\\
  \hline
  $b_{12}$ & $N(50,10^2)$ & $N(300,10^2)$ \\
  $b_{12}$ & $N(300,10^2)$ & $N(300,100^2)$ \\
  \hline
\end{tabular}
\end{center}
\end{table}

\begin{proposition} (Observation 4)
A candidate node without historical demand records can be chosen to locate an airport.
\end{proposition}

Node 3 is chosen to locate an additional airport when we increase the upper bound of simulated demand for O-D pair between node 1 and node 3.
This observation is intuitive, because the size of confidence set increases as the upper bound increases, leading to a new worst-case distribution. A direct implication is to decide whether to setup airport close to potential future customers who have never requested logistics service before.

\begin{figure}[H]
 \centering
 \includegraphics[width=0.80\textwidth]{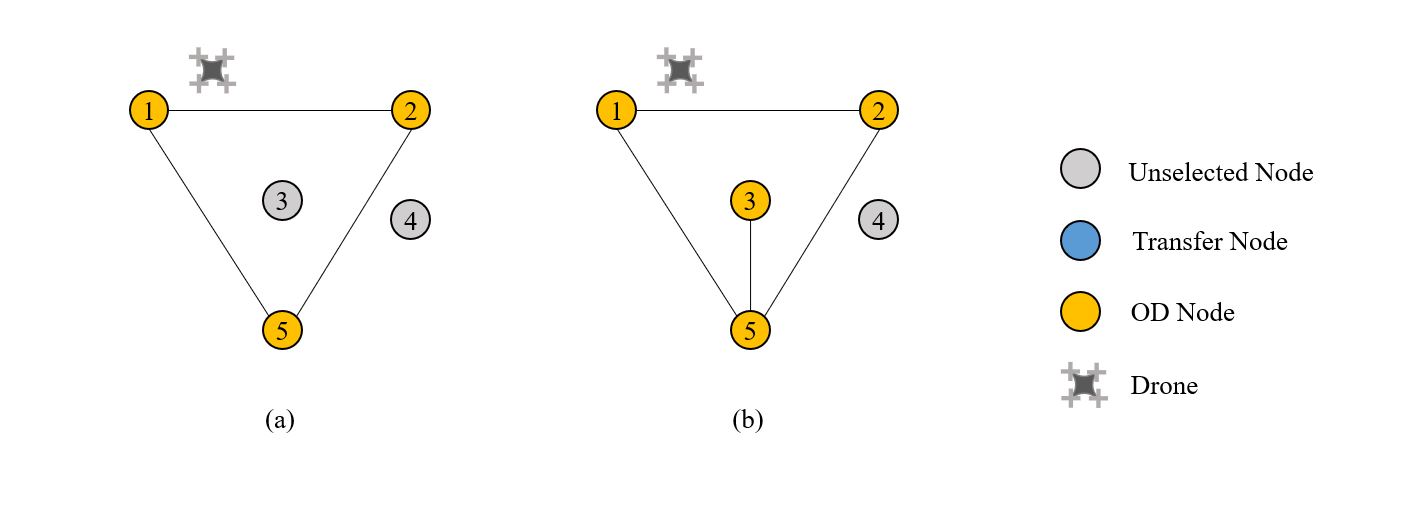}
 \caption{Demand occurs with zero demand}
\end{figure}

We summarize our results in the following observations:
\begin{itemize}
  \item The marginal decision of setting up an additional airport at a particular candidate location depends crucially on both the infrastructure cost at and transportation cost en route to the location.
  \item The optimal network configuration depends jointly on infrastructural and transportation cost coefficients,  the battery capacity constraints, as well as the pooling effects in channel capacities.
  \item The optimal network design is relatively robust under demand uncertainty.
  \item A candidate node without historical demand records can be chosen to locate an airport.
\end{itemize}

\subsection{Case study}

We perform a case study of UAMN design problem wherein the historical demand data is supplied by our industry partner and the parameters are validated accordingly. We set $\beta=50$, $\theta=100$, $L=18$, $t=1$, $y$ is binary.
Transportation cost coefficients and battery consumption are proportional to their corresponding distance.
Currently, airports are standardized products and we assume $Cs$ is homogeneous. We also set $Cf$ and $Cd$ identical across different locations in this study.

Table 4 provides historical demand data. There are 7 blood collection points in Hangzhou, China, including Wushan Square, Longxiang Mansion, Hangzhou Theater, Wulin Courtyard, Zheyi Blood Square, Yunhe Square, and Xiasha Wu Mart. All these 7 blood collection points supply blood to blood bank for processing every day. We assume these 7 O-D pairs all take the upper bound of demand as 25, and take the lower bound od demand as 0. We assume all other O-D pairs take both the upper and the lower bound of demand as zeros in this pilot study.

\begin{table}[H]
\caption{Blood Transportation Volume (kg) Depend on Location and Date}
\begin{center}
\resizebox{\textwidth}{28mm}{
\begin{tabular}{c c c c c c c c}
  \hline
  Date & 2018/11/12 & 2018/11/13 &2018/11/14 &2018/11/15 &2018/11/16  &2018/11/17 &2018/11/18 \\
  \hline
  Wushan Square        & 0.4 & 3.1 & 4.6 & 3.1 & 1.3 & 6.4 & 1.5\\
  Longxiang Mansion    & 4.3 & 4.2 & 4.4 & 4.93& 5.2 & 9.9 & 8.4\\
  Hangzhou Theater     & 1.3 & 1.2 & 1.4 & 0   & 0.8 & 3.2 & 0.5\\
  Wulin Courtyard      & 5.2 & 9.9 & 5.5 & 5.9 & 0.5 & 12.3&12.7\\
  Zheyi Blood Station & 1 & 7 & 4.2 & 4.5 & 10.7 & 0 & 0\\
  Yunhe Square & 0 & 4.8 & 0 & 5.1 & 0 & 5.4 & 3.6\\
  Xiasha Wu Mart & 0 & 5.2 & 0 & 2.4 & 0 & 4.82 & 5.5\\
  \hline
\end{tabular}}
\end{center}
\end{table}

Distance data between all blood collection points and blood bank are given in the Appendix. The main problem for the company is that Xiasha Wu Mart is far away from the blood bank, which outlasts the battery capacity of their drones. Two options are to be evaluated: (A) Starting from this remote location, the drone make a detour and transit at another location for charging. This option incurs higher transportation costs due to the longer flight. We also need to decide the transition point with sufficient capacity. (B) We set up an additional candidate airport, which can potentially save the transportation cost but require additional infrastructure investment.

We solve this problem and attain the worst-case distribution of demand which are displayed in Appendix. We can notice that the worst case distribution of demand is different from the original sample demand. The value of demand can be their upper bound, lower bound, or themselves. This table illustrates that our model considers the worst case under the confidence lever of $\beta$.
Figure 4 shows the blood transportation network design. Xiasha Wu mart and blood center are linked by Zheyi blood donation station, rather than the candidate point Shenbo paradise, which means we eventually choose option (A) or (B).

\begin{figure}[H]
 \centering
 \includegraphics[width=0.80\textwidth]{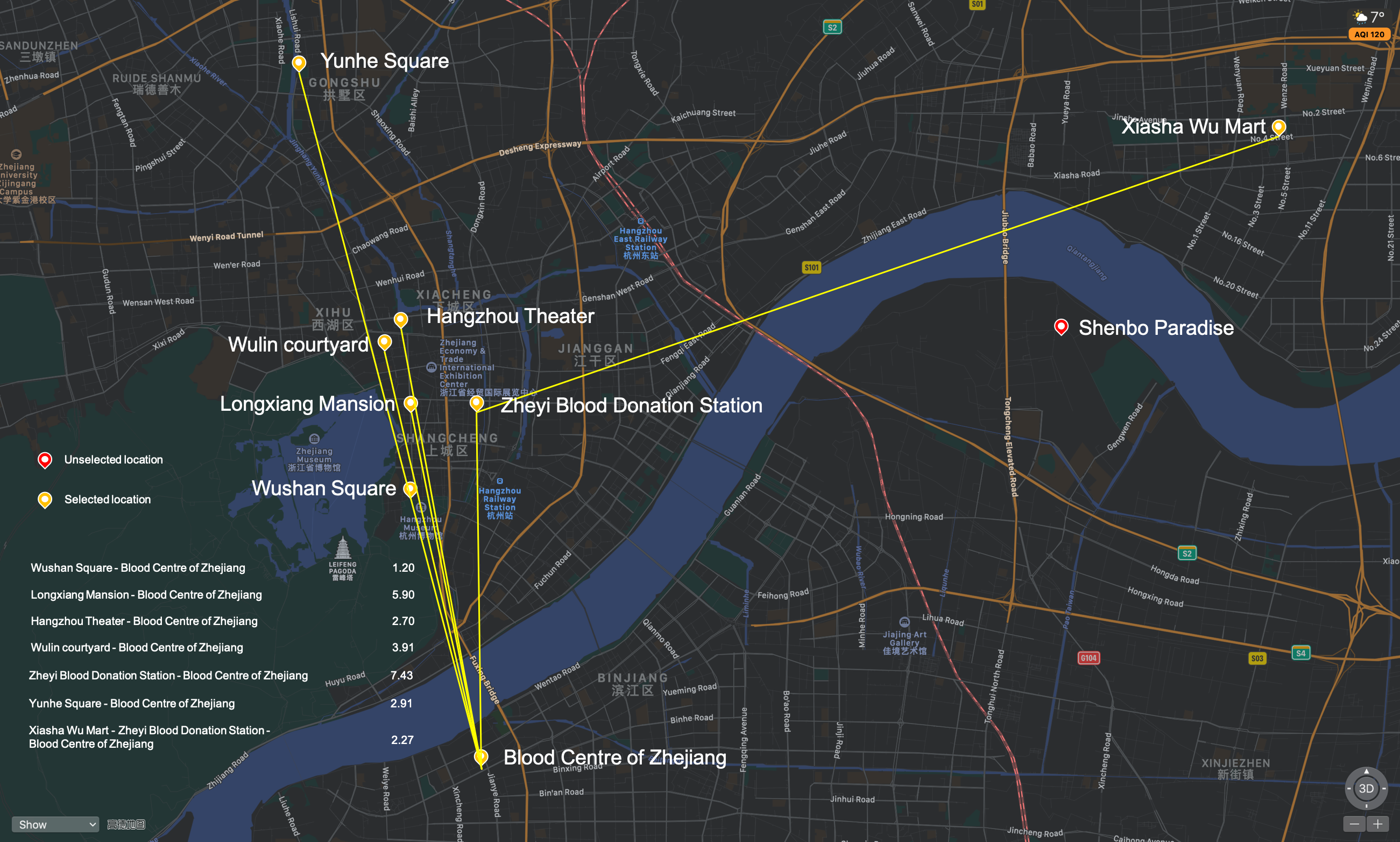}
 \caption{Blood transportation network}
\end{figure}

In summary, our model framework provides an exact approach to design a UAMN for drone delivering. We can easily apply it on other real case studies related to drone delivering beyond blood transportation.

\section{Conclusion}

\label{s-con}
In this paper, we propose an integrated facility location and network design problem, i.e., UAMN design problem for drone logistics. To do this, we develop a risk-averse two-stage stochastic model that is based on Wasserstein distance. We then develop a reformulation technique that simplifies the worst-case expectation term in the original model, and obtain a Min-Max model correspondingly. By using Lagrange multipliers, we successfully decompose decision variables and reduce the computational complexity, allowing us to solve real-world instances under our model framework.
A few managerial observations are obtained through numerical examples. For example, we find that the optimal network configuration is affected by the ``pooling effects" in channel capacities. A nice feature of our DRO framework is that the optimal network design is relatively robust under demand uncertainty. Interestingly, a candidate node without historical demand records can be chosen to locate an airport. We next demonstrate the application of our model by providing a real problem involving medical resources transportation with our industry partner.


\bibliographystyle{chicago}
\bibliography{RML}

\newpage
\appendix
\noindent\textbf{\Large{Online Appendices}}
\section{Proof of Propositions}
\subsection{Proof for Proposition 1}
\begin{proof}
We cite the lemma in \cite{zhao2018data} as follows:
\begin{lemma}
Assuming that there are N historical data samples $\xi^1, \xi^2, \cdots, \xi^N$ which are i.i.d drawn from the true distribution P, for any fixed first-stage decision x, we have
\begin{equation*}
 \max_{\widehat{P}\in \mathcal{D}}E_{\widehat{P}}[\mathcal{Q}(x,\xi)]=\min_{\beta\geq0}\frac{1}{N}\sum_{i=1}^{N}\max_{\xi}[{\mathcal{Q}(x,\xi)-\beta\rho(\xi,\xi^i)}]+\theta\beta
\end{equation*}
\end{lemma}
where $x$ is the decision variables, and $\xi$ is the sample which is uncertainty. Our decision variables are $\eta$, and the demand $b_k$ is uncertainty. We substitute $x$ as $\eta$, and substitute $\xi$ as $b_k$. Then, the conclusion of proposition 1 is attained.
\end{proof}

\subsection{Proof for Proposition 2}
\begin{proof}
We first take the dual of the formulation for the second-stage cost (i.e., $\mathcal{Q}(\eta,b)$ and combine it with the second problem to obtain the following subproblem (denoted as SUB) corresponding to each sample $b^j,j=1,\cdots,N$:

   \begin{align*}
        (SUB)\phi^j(\eta,\beta)&=\max_{b}\quad \{\mathcal{Q}(\eta,b)-\beta\rho^{j}(b)\}\\
        &=\max_{b,\mu,\lambda\geq0}[\min_{x}\{(C^T+\mu^TA+\lambda^TD)X-\mu^TB(b)-\lambda^TE\}-\beta\rho^{j}(b)]\\
        &=\max_{b,\mu,\lambda\geq0}[-\mu^TB(b)-\lambda^TE-\beta\rho^{j}(b)]\\
        \text{s.t.}\quad &C^T+\mu^TA+\lambda^TD\geq0
   \end{align*}
 We argue that $C^T+\mu^TA+\lambda^TD\geq0$. If argument does not hold, let $X$ tend to positive infinity, then $(C^T+\mu^TA+\lambda^TD)X$ tends to negative infinity. Because we want to find max over $\mu$, and$\lambda$, $C^T+\mu^TA+\lambda^TD\geq0$ holds.
The third equation holds because, if $(C^T+\mu^TA+\lambda^TD)_i=0$, then $x_i$ can be any nonnegative number, if $(C^T+\mu^TA+\lambda^TD)_i>0$, then $x_i$ will be zero, so $(C^T+\mu^TA+\lambda^TD)X=0$ always holds when we min over $X$.

Let $\rho^{j}(b)=\sum_{i=1}^{K}|b_i-b_i^j|$, for a fixed $\lambda$,$\mu$, obtain an optimal solution $b$ to the problem (SUB):
\begin{equation*}
\phi^j(\eta,\beta,\lambda,\mu)=\max_{b}[-\sum_{i=1}^{K}\mu_{i}b_i+\sum_{i=1}^{K}\mu_{K+i}b_i-\beta\sum_{i=1}^K|b_i-b_i^j|-\lambda^TE]
\end{equation*}
So it can be observed that at least one optimal solution $b^*$ to the subproblem (SUB) satisfies
$b_i^*=W_i^{-}$,$b_i^*=W_i^{+}$,or $b_i^*=b_i^j$, for each $i=1,2,\cdots,K$,
which indicate that for each $i=1,\cdots,K$, the $i$th component of optimal solution $b^*$ can achieve its lower bound $W_i^{-}(\delta_i^{+}=0,\delta_i^{-}=1)$, upper bound $W_i^{+}(\delta_i^{+}=1,\delta_i^{-}=0)$, or the sample value $b_i^j(\delta_i^{+}=\delta_i^{-}=0)$.

Then, there exists an optimal solution $b^*$ of (SUB) satisfying the following constraints:

\begin{equation*}
b_i=(W_i^{+}-b_i^j)\delta_i^{+}+(W_i^{-}-b_i^j)\delta_i^{-}+b_i^j, \forall i=1,\cdots,m,
\end{equation*}
\begin{equation*}
\delta_i^{+}+\delta_i^{-}\leq1, \delta_i^{+},\delta_i^{-}\in\{0,1\}, \forall i=1,\cdots,m.
\end{equation*}

We can linearize the $\phi^j(\eta,\lambda,\mu)$ as:

\begin{align*}
\phi^j(\eta,\beta,\lambda,\mu)=&\max_{b}[-\sum_{i=1}^{K}(\mu_{i}-\mu_{K+i})b_i-\beta\sum_{i=1}^K|b_i-b_i^j|-\lambda^TE]\\
=&\max_{b}[-\sum_{i=1}^{K}(\mu_{i}-\mu_{K+i})((W_i^{+}-b_i^j)\delta_i^{+}+(W_i^{-}-b_i^j)\delta_i^{-}+b_i^j)\\
&-\beta\sum_{i=1}^K|(W_i^{+}-b_i^j)\delta_i^{+}+(W_i^{-}-b_i^j)\delta_i^{-}|-\lambda^TE]\\
=&\max_{b}[-\sum_{i=1}^{K}(\mu_{i}-\mu_{K+i})((W_i^{+}-b_i^j)\delta_i^{+}+(W_i^{-}-b_i^j)\delta_i^{-}+b_i^j)\\
&-\beta\sum_{i=1}^K((W_i^{+}-b_i^j)\delta_i^{+}-(W_i^{-}-b_i^j)\delta_i^{-})-\lambda^TE]
\end{align*}

The (SUB) can be described as follows:
\begin{align*}
\varpi^j(\eta,\beta)
=&\max_{\delta^{+},\delta^{-},\mu,\lambda\geq0}[-\sum_{i=1}^{K}(\mu_{i}-\mu_{K+i})((W_i^{+}-b_i^j)\delta_i^{+}+(W_i^{-}-b_i^j)\delta_i^{-}+b_i^j)\\
&-\beta\sum_{i=1}^K((W_i^{+}-b_i^j)\delta_i^{+}-(W_i^{-}-b_i^j)\delta_i^{-})-\lambda^TE]\\
\text{s.t.}\quad &C^T+\mu^TA+\lambda^TD\geq0\\
&\delta_i^{+},\delta_i^{-}\in\{0,1\}\\
&\delta_i^{+}+\delta_i^{-}\leq1.
\end{align*}
\end{proof}

\section{Matrix Definitions}
The second-stage problem $\mathcal{Q}(\eta,b)$ is:
   \begin{align*}
        \min_{X}\quad &C^{T}X\\
        \text{s.t.}\quad &AX=B\\
        &DX\leq E,
   \end{align*}
\begin{itemize}
  \item $B=(b^T,-b^T,0^T)^T$,$b=(b_1,b_2,\cdots,b_K)^T$. $B$ is a vector, and its size is determined by $K$.
  \item $A$ is a matrix composed by coefficient before $X$ in constraints (1),(2),(3). The number of row is determined by the number of constraints (1) (2) (3), and the number of column is determined by the size of vector $X$, i.e. it is equal to the number of elements in vector $X$.
  \item $D$ is a matrix composed by coefficient before $X$ in constraints (4),(5),(6). The number of row is determined by the number of constraints (4) (5) (6), and the number of column is determined by the size of vector $X$, i.e. it is equal to the number of elements in vector $X$.
  \item $E$ is a vector composed by right side in constraints (4),(5),(6). The number of row is determined by the number of constraints (1) (2) (3).
  \item $C$ is a vector composed by coefficients for $X$ in objective function. The number of row is determined by the number of elements in vector $X$.
\end{itemize}

\section{Data for Case Study}
We give the pairwise distances between any two candidate locations in Table 5.
\begin{table}[H]
\caption{Pairwise Distance Between two Candidate Locations}
\begin{center}
\resizebox{\textwidth}{28mm}{
\begin{tabular}{c c c c c c c c c c}
  \hline
  Distance (km) & Wushan Square  & Longxiang Mansion & Hangzhou Theater &Wulin Courtyard  &Zheyi Blood Station &Yunhe Square & Xiasha Wu Mart & Blood Center & Candidate Point\\
  \hline
  Wushan Square        & 0   & 1.8 & 3.4 & 3   & 2.1 & 8.9 & 19   & 5.3 & 13.5\\
  Longxiang Mansion    & 1.8 & 0   & 1.5 & 1.2 & 1.3 & 7.1 & 18.4 & 7   & 13.1\\
  Hangzhou Theater     & 3.4 & 1.5 & 0   & 0.5 & 2.2 & 5.6 & 18.2 & 8.6 & 13.3\\
  Wulin Courtyard      & 3   & 1.2 & 0.5 & 0   & 2.1 & 5.9 & 18.5 & 8.2 & 13.5\\
  Zheyi Blood Station  & 2.1 & 1.3 & 2.2 & 2.1 & 0   & 7.8 & 17.1 & 6.7 & 11.8\\
  Yunhe Square         & 8.9 & 7.1 & 5.6 & 5.9 & 7.8 & 0   & 19.9 & 14.2& 16.3\\
  Xiasha Wu Mart       & 19  & 18.4& 18.2& 18.5& 17.1& 19.9& 0    & 20.3& 5.9\\
  Blood Center         & 5.3 & 7   & 8.6 & 8.2 & 6.7 & 14.2& 20.3 & 0   & 14.4\\
  Candidate Point      & 13.5& 13.1& 13.3& 13.5& 11.8& 16.3& 5.9  & 14.4& 0\\
  \hline
\end{tabular}}
\end{center}
\end{table}

The worst-case demand distributions are given in Table 6.
\begin{table}[H]
\caption{Worst demand (kg) distribution}
\begin{center}
\resizebox{\textwidth}{28mm}{
\begin{tabular}{c c c c c c c c}
  \hline
  Date & 2018/11/12 & 2018/11/13 &2018/11/14 &2018/11/15 &2018/11/16  &2018/11/17 &2018/11/18 \\
  \hline
  Wushan Square        & 0.4 & 3.1 & 4.6 & 3.1 & 1.3 & 25 & 1.5\\
  Longxiang Mansion    & 4.3 & 4.2 & 4.4 & 4.93& 5.2 & 25 & 8.4\\
  Hangzhou Theater     & 1.3 & 1.2 & 1.4 & 0   & 0.8 & 25 & 0.5\\
  Wulin Courtyard      & 5.2 & 9.9 & 5.5 & 5.9 & 0.5 & 25 &12.7\\
  Zheyi Blood Station & 1 & 7 & 4.2 & 4.5 & 25 & 0 & 0\\
  Yunhe Square & 0 & 4.8 & 0 & 5.1 & 0 & 25 & 3.6\\
  Xiasha Wu Mart & 0 & 5.2 & 0 & 2.4 & 25 & 25 & 5.5\\
  \hline
\end{tabular}}
\end{center}
\end{table}

\section{Derivation for Computation Technique}
In this section, we provide further explanations why our solution procedure works. Let
\begin{align*}
\zeta^j=&-\sum_{i=1}^{K}(\mu_{i}^j-\mu_{K+i}^j)((W_i^{+}-b_i^j)\delta_i^{j+}+(W_i^{-}-b_i^j)\delta_i^{j-}+b_i^j)\\
&-\beta\sum_{i=1}^K((W_i^{+}-b_i^j)\delta_i^{j+}-(W_i^{-}-b_i^j)\delta_i^{j-})-(\lambda^j)^TE
\end{align*}%

Then change the order of max and min:
\begin{align*}
\psi=\min_{\eta,\beta\geq0}\max_{\delta^{+},\delta^-,\mu,\lambda\geq0}\quad&\frac{1}{N}\sum_{j=1}^{N}\zeta^j+\theta\beta+\sum_{t\in
T}\sum_{(i,j)\in A}Cd_{ij}^{t}y_{ij}^{t}+\sum_{j\in V}Cf_{j}z_{j}+\sum_{j\in
V}Cs_{j}w_{j}z_{j}\\
=\max_{\delta^{+},\delta^-,\mu,\lambda\geq0}\min_{\eta,\beta\geq0}\quad &\frac{1}{N}\sum_{j=1}^{N}\zeta^j+\theta\beta+\sum_{t\in
T}\sum_{(i,j)\in A}Cd_{ij}^{t}y_{ij}^{t}+\sum_{j\in V}Cf_{j}z_{j}+\sum_{j\in
V}Cs_{j}w_{j}z_{j}\\
\text{s.t.}\quad  &C^T+(\mu^j)^TA+(\lambda^j)^TD\geq0\\
&\delta_i^{+},\delta_i^{-}\in\{0,1\}\\
&\delta_i^{+}+\delta_i^{-}\leq1.
\end{align*}%
where $\eta=(z,y)$.

Separate the matrix $E$, and notice that $E$ contains the right side of (4), (5), (6). We separate $\frac{1}{N}\sum_{j=1}^{N}-(\lambda^j)^TE$:
\begin{align}
\psi=&\max_{\delta^{+},\delta^-,\mu,\lambda\geq0}\min_{\eta,\beta\geq0}\quad -\frac{1}{N}\sum_{j=1}^{N}
[\sum_{i=1}^{K}(\mu_{i}^j-\mu_{K+i}^j)((W_i^{+}-b_i^j)\delta_i^{j+}+(W_i^{-}-b_i^j)\delta_i^{j-}+b_i^j)]\notag\\
+&[-\frac{1}{N}\sum_{j=1}^{N}
\sum_{i=1}^K((W_i^{+}-b_i^j)\delta_i^{j+}-(W_i^{-}-b_i^j)\delta_i^{j-})+\theta]\beta\notag\\
+&\sum_{j\in V}[(-(\overline{\lambda}_{1})_j(w_j-M)+Cf_{j}+Cs_{j}w_{j})z_j-(\overline{\lambda}_{1})_jM]\notag\\
+&\sum_{t\in T}\sum_{(i,j)\in A}(Cd_{ij}^{t}-(\overline{\lambda}_{2})_{ij}u_{ij}^t)y_{ij}^{t}\notag\\
+&\sum_{j=1}^{N}(C^T+(\mu^j)^TA+(\lambda^j)^TD)\gamma^j\notag\\
\text{s.t.}\quad
&\delta_i^{+},\delta_i^{-}\in\{0,1\}\notag\\
&\delta_i^{+}+\delta_i^{-}\leq1.\notag
\end{align}
where $\lambda_1$ represents the multiplier of constraints (5), $\lambda_2$ represents the multiplier of constraints (4), and $\lambda_3$ represents the multiplier of constraints (6).  So the multiplier $\lambda$ can be separated by $\lambda=({\lambda_2}^T,{\lambda_1}^T,{\lambda_3}^T)^T$. $\lambda_1^j$ is corresponding to $j$th sample of $\lambda_1$, $\lambda_2^j$ is corresponding to $j$th sample of $\lambda_2$, and $\lambda_3^j$ is corresponding to $j$th sample of $\lambda_3$.
 $\overline{\lambda}_1$ is the mean of $\lambda_1^j$  , and $\overline{\lambda}_2$ is the mean of $\lambda_2^j$. So the multiplier $\lambda^j$ can be separated by $\lambda^j=({\lambda_2^j}^T,{\lambda_1^j}^T,{\lambda_3^j}^T)^T$.
$(\overline{\lambda}_{1})_j$ is the component in $\overline{\lambda}_{1}$ corresponding to $z_j$. $(\overline{\lambda}_{2})_{ij}$ is the component in $\overline{\lambda}_{2}$ corresponding to $y_{ij}$.

Separate the last term in above formula:
\begin{align*}
\psi
=&\max_{\delta^{+},\delta^-,\mu}\min_{\beta\geq0}\quad \{-\frac{1}{N}\sum_{j=1}^{N}
[\sum_{i=1}^{K}(\mu_{i}^j-\mu_{K+i}^j)((W_i^{+}-b_i^j)\delta_i^{j+}+(W_i^{-}-b_i^j)\delta_i^{j-}+b_i^j)-(\mu^j)^TA\gamma^j]\\
+&[-\frac{1}{N}\sum_{j=1}^{N}
\sum_{i=1}^K((W_i^{+}-b_i^j)\delta_i^{j+}-(W_i^{-}-b_i^j)\delta_i^{j-})+\theta]\beta\}\\
+&\max_{\lambda\geq0}\min_{\eta}\{
\sum_{j\in V}[(-(\overline{\lambda}_{1})_j(w_j-M)+Cf_{j}+Cs_{j}w_{j})z_j-(\overline{\lambda}_{1})_jM]+\sum_{j=1}^{N}(\lambda_1^j)^TD_1\gamma^j
+\sum_{j=1}^{N}(\lambda_3^j)^TD_3\gamma^j\}\\
+&\max_{\lambda\geq0}\min_{\eta}\{
\sum_{t\in T}\sum_{(i,j)\in A}(Cd_{ij}^{t}-(\overline{\lambda}_{2})_{ij}u_{ij}^t)y_{ij}^{t}+\sum_{j=1}^{N}(\lambda_2^j)^TD_2\gamma^j\}\\
+&\sum_{j=1}^{N}C^T\gamma^j\\
\text{s.t.}\quad
&\delta_i^{+},\delta_i^{-}\in\{0,1\}\\
&\delta_i^{+}+\delta_i^{-}\leq1.
\end{align*}
where $D_1$ represents the coefficients in constraints (5), $D_2$ represents the coefficients in constraints (4), and $D_3$ represents the coefficients in constraints (6).  So the multiplier $D$ can be separated by $D=({D_2}^T,{D_1}^T,{D_3}^T)^T$. $(D_{2})_{(ij)}$ is a row of $D_2$ corresponding to $y_{ij}$.

We can easily find that every problem of Max-Min is independent because the parameters in each problem are different. After we get the results, we need to check if the constraint holds. To see this, we first change the order of sum and max-min:

\begin{align*}
\psi
=&\max_{\delta^{+},\delta^-,\mu}\min_{\beta\geq0}\{-\frac{1}{N}\sum_{j=1}^{N}
[\sum_{i=1}^{K}(\mu_{i}^j-\mu_{K+i}^j)((W_i^{+}-b_i^j)\delta_i^{j+}+(W_i^{-}-b_i^j)\delta_i^{j-}+b_i^j)-(\mu^j)^TA\gamma^j]\\
+&[-\frac{1}{N}\sum_{j=1}^{N}
\sum_{i=1}^K((W_i^{+}-b_i^j)\delta_i^{j+}-(W_i^{-}-b_i^j)\delta_i^{j-})+\theta]\beta\}\tag{a}\\
+&\max_{\lambda_1,\lambda_3\geq0}\min_{z}\{
\sum_{j\in V}[(-(\overline{\lambda}_{1})_j(w_j-M)+Cf_{j}+Cs_{j}w_{j})z_j-(\overline{\lambda}_{1})_jM]+\sum_{j=1}^{N}(\lambda_1^j)^TD_1\gamma^j
+\sum_{j=1}^{N}(\lambda_3^j)^TD_3\gamma^j\}\tag{b}\\
+&\sum_{(i,j)\in A}\max_{\lambda_{2ij}\geq0}\min_{y_{ij}^t}\{
[\sum_{t\in T}(Cd_{ij}^{t}-(\overline{\lambda}_{2})_{ij}u_{ij}^t)y_{ij}^{t}+\sum_{h=1}^{N}(\lambda_{2ij}^h)^TD_{2ij.}\gamma^h]\}\tag{c}\\
+&\sum_{j=1}^{N}C^T\gamma^j\\
\text{s.t.}\quad
&\delta_i^{+},\delta_i^{-}\in\{0,1\}\\
&\delta_i^{+}+\delta_i^{-}\leq1.
\end{align*}
We must remember that the multiplier $\gamma$ represents the condition $C^T+(\mu^j)^TA+(\lambda^j)^TD\geq0$.
So we want to find a proper  $\gamma$ to make this constraint hold. To summarize, we give a description of how to solve every Max-Min problem. \\
\begin{itemize}
  \item Given $\beta$, (a) is easy to solve because it is a linear Max problem .
  \item (b) cannot be separated into independent terms only related to one $z_j$, because there are $z$ in $D_3$. The number of $z$ is usually not too large, so solving the (b) is acceptable. Given the value of $z$, it is a linear Max problem.

  \item Given $y_{ij}^{t}$, it is a linear Max problem. If $y_{ij}^{t}$ is a binary variable, give a example of (c).

Q represents
\begin{align*}
\max_{\lambda_{2ij}\geq0}\min_{y_{ij}^t}\{
[\sum_{t\in T}(Cd_{ij}^{t}-\overline{\lambda}_{2ij}u_{ij}^t)y_{ij}^{t}+\sum_{h=1}^{N}(\lambda_{2ij}^h)^TD_{2ij.}\gamma^h]\}
\end{align*}
If $(Cd_{ij}^{t}-\overline{\lambda}_{2ij}u_{ij}^t)<0$, then $y_{ij}^{t}=1$;
If $(Cd_{ij}^{t}-\overline{\lambda}_{2ij}u_{ij}^t)\geq0$, then $y_{ij}^{t}=0$.

\end{itemize}

\end{document}